\documentclass[11pt]{amsart}
\usepackage{amssymb,amsmath,amsfonts,amscd,euscript}
\usepackage{tikz,graphicx}
\usepackage{tikz-cd}
\usepackage[all]{xy}
\usepackage{color}

\usepackage[backref=page]{hyperref}
\usepackage{todonotes}

\usepackage{amsthm,enumitem,bm,pict2e,xcolor}
\usepackage{indentfirst}

\usepackage{algorithm}
\usepackage{algpseudocode}
\usepackage{float}

\newtheorem{theoremA}{Theorem}

\newtheorem{conjA}[theoremA]{Conjecture}

\newcommand{\nc}{\newcommand}

\numberwithin{equation}{section}
\newtheorem{thm}{Theorem}[section]
\newtheorem{prop}[thm]{Proposition}

\newtheorem{rem}[thm]{Remark}

\newtheorem{example}[thm]{Example}

\newtheorem{conj}[thm]{Conjecture}

\nc{\gl}{\mathfrak{gl}}
\nc{\GL}{\mathfrak{GL}}
\nc{\g}{\mathfrak{g}}
\nc{\gh}{\widehat\g}
\nc{\hh}{{\fh}^a}
\nc{\fh}{\mathfrak{h}}
\nc{\la}{\lambda}
\nc{\al}{\alpha }
\nc{\be}{\beta }
\nc{\om}{\omega }

\nc{\ta}{\theta}
\nc{\veps}{\varepsilon}
\nc{\ch}{{\mathop {\rm ch}}}
\nc{\Tr}{{\mathop {\rm Tr}\,}}
\nc{\Id}{{\mathop {\rm Id}}}
\nc{\Aut}{{\mathop {\rm Aut}}}
\nc{\End}{{\mathop {\rm End}}}

\nc{\bra}{\langle}
\nc{\ket}{\rangle}
\nc{\bs}{{\bf s}}
\nc{\bp}{{\bf p}}

\nc{\pa}{\partial}
\nc{\ld}{\ldots}
\nc{\cd}{\cdots}
\nc{\hk}{\hookrightarrow}
\nc{\T}{\otimes}
\nc{\Gr}{\mathrm{Gr}}
\nc{\aGr}{{\mathcal Gr}}
\nc{\aFl}{{\mathcal Fl}}
\nc{\ov}{\overline}

\nc{\msl}{\mathfrak{sl}}
\nc{\mgl}{\mathfrak{gl}}
\nc{\U}{\mathrm U}
\nc{\Res}{\mathrm{Res\ }}

\newcommand{\bC}{{\mathbb C}}

\newcommand{\bZ}{{\mathbb Z}}

\newcommand{\bP}{{\mathbb P}}

\newcommand{\eO}{\EuScript{O}}

\newcommand{\Hom}{\mathrm{Hom}}
\newcommand{\ba}{{\bf a}}
\newcommand{\bfm}{{\bf m}}

\newcommand{\bd}{{\bf d}}
\newcommand{\bk}{{\mathbf k}}

\tikzcdset{scale cd/.style={every label/.append style={scale=#1},
		cells={nodes={scale=#1}}}}
\begin{document}

\title[PrIncipal quiver Grassmannians: conjectures]
{PrIncipal quiver Grassmannians: conjectures}

\author{Stanislav Fedotov}
\address{Nebius, London, UK}
\email{st.n.fedotov@gmail.com}

\author{Evgeny Feigin}
\address{School of Mathematical Sciences, Tel Aviv University, Tel Aviv, 69978, Israel}
\email{evgfeig@gmail.com}

\begin{abstract}
Let $P$ and $I$ be a projective and an injective representations of a Dynkin quiver.
We consider quiver Grassmannians of subrepresentations of dimension $\dim P$ inside
representations of dimension $\dim P + \dim I$. 	
Based on extensive computer experiments, we formulate several conjectures about the 
algebro-geometric properties of these quiver Grassmannians. 	
\end{abstract}

\maketitle

\section{Introduction}
Let $Q$ be a Dynkin quiver with an arbitrary orientation. Let $P$ and $I$ be a projective
and an injective representations of $Q$. Our main objects of study are the quiver Grassmannians
$\Gr_{\dim P}(M)$, where $M$ are $\dim P + \dim I$ dimensional representations of $Q$. 

Examples of such varieties naturally show up in various problems of algebraic geometry and representation theory 
\cite{CI20,CFFFR17,CFFFR20,CFR12}.
In particular, classical flag  varieties and their PBW degenerations can be realized in this way (see \cite{CL15,Fe11,Fe12-1,FeFi13,FFR17}).
Our goal
is twofold: first, we develop the software which allows to perform calculations with the 
quiver representations and quiver Grassmannians. Second, we use the software to check various
conjectural properties of our quiver Grassmannians. In what follows, following the suggestion 
of Markus Reineke, we call quiver Grassmannians $\Gr_{\dim P}(M)$ with $\dim M =\dim P +\dim I$ principal (or PrIncipal), where 
Pr stands for projective and In for injective (we note that
one can not study all possible quiver Grassmannians, since every projective variety can be realized in this way \cite{Re13,Ri18}).

Let $Q_0$ and $Q_1$ be the sets of vertices and arrows of $Q$. Let $\bd=\dim P+\dim I$ 
be the dimension vector of a representations $M$ as above and  
let ${\rm Rep}_\bd$ be the representation space of $Q$ of dimension $\bd$. A point of the 
universal quiver Grassmannian $\Gr_{\dim P}(\bd)$ consists of a point 
$M\in  {\rm Rep}_\bd$ and a collection of $(\dim P)_i$-dimensional subspaces of $M_i$, $i\in Q_0$
compatible with the maps of $M$. One has a natural projection map
$\pi: \Gr_{\dim P}(\bd)\to {\rm Rep}_\bd$ whose fibers are the  quiver Grassmannians.  
Our first goal is to describe the reduced scheme structure of these fibers.

By definition, $\Gr_{\dim P}(M)$ are embedded into the product of classical Grassmann varieties 
$\Gr_{\dim P_i}(M_i)$, $i\in Q_0$ and hence admit the Pl\"ucker embedding into the product of
projectivized wedge powers. 
One knows \cite{LW19} that the quiver Grassmannians can be endowed with the scheme structure defined
by quadratic relations labeled by the edges $\alpha\in Q_1$. However, in most cases this scheme structure
is not reduced. One can also naturally define quadratic relations for each path in $Q$. Here is our first
conjecture.

\begin{conjA}\label{conj:redscheme}
For every principal quiver Grassmannian the scheme structure defined by the quadratic relations 
corresponding to all paths in $Q$ is reduced. 
\end{conjA}  

The representation space ${\rm Rep}_\bd$ is naturally acted upon by the group $G_\bd$, the product
of $GL_{d_i}$ for all $i\in Q_0$. Since $Q$ is Dynkin, there is a unique open orbit for this action;
let $M^0$ be an element of this open orbit. 
Then $\Gr_{\dim P}(M^0)$ is irreducible and of dimension
$\langle \dim P,\dim I\rangle$; hence for any $\bd$-dimensional $M$ the dimension of $\Gr_{\dim P}(M^0)$ 
is at least $\langle \dim P,\dim I\rangle$
(see \cite{CFR12}). We say that $M$ degenerates to $N$ ($M,N\in {\rm Rep}_\bd$),
if $N$ is contained in the closure of the $G_\bd$ orbit of $M$.   
 We put forward the following conjecture.

\begin{conjA}\label{conj:M^1}
There exists a representation 
$M^1$ of dimension $\dim P$ such that a quiver Grassmannian 
$\Gr_{\dim P}(N)$, $N\in {\rm Rep}_\bd$ is irreducible of dimension $\langle \dim P,\dim I\rangle$ if and only if $N$ degenerates to $M^1$.
If $\dim P$ and $\dim Q$ have no zero components, then $M^1 = P\oplus I$.
\end{conjA}  
   
From \cite{CFR13-1}, we already know that any representation degenerating to $P\oplus I$ is of dimension $\langle \dim P,\dim I\rangle$. Conjecture \ref{conj:M^1} is proved in \cite{CFFFR17,CFFFR20} for equioriented type $A$ quiver in certain special cases.

Another natural question to ask is for which representations $N$
the quiver Grassmannians $\Gr_{\dim P}(N)$ are of the expected (minimal possible)
dimension $\langle \dim P,\dim I\rangle$. 
Some special cases of this question were investigated in \cite{CFFFR17}, \cite{CFFFR20}. Here is our third conjecture.

\begin{conjA}\label{conj:M^2}
Let $Q$ be of type $A$ with arbitrary orientation. There exists a representation 
$M^2$ of dimension $\dim P$ such that a quiver Grassmannian 
$\Gr_{\dim P}(N)$, $N\in {\rm Rep}_\bd$ is of dimension $\langle \dim P,\dim I\rangle$ if and only if $N$ degenerates to $M^2$.
\end{conjA}  
 
In words, Conjecture \ref{conj:M^2} states that there exists the deepest quiver 
Grassmannian of the minimal possible dimension. We give a conjectural description of $M^2$ in some special cases.
We note that if Conjecture \ref{conj:M^2} holds true then the subvariety of the 
universal quiver Grassmannian consisting of fibers over representations degenerating
to $M^2$ is flat (\cite{CFFFR17,CFFFR20}).
We also note that Conjecture \ref{conj:M^2} does not hold in general: in type $D_4$ for certain
representations $P$ and $I$ there exist 
three deepest quiver Grassmannians of the minimal possible dimension.

A universal description of minimal dimension locus is conjecturally given, under quite lax restrictions, in terms of dimensions of homomorphism spaces. Here's our next conjecture:

\begin{conjA}\label{conj:hom-dim}
Let either of the following two conditions holds:
a)\ $Q$ is of type $A$ and neither $\dim P$ nor $\dim I$ have zero components,
b)\ $Q$ is of type $D$ and $P\oplus I$ contains  every indecomposable projective and every indecomposable injective representation as a summand.
Then $\dim\Gr_{\dim P}(M)=\langle \dim P,\dim I\rangle$ if and only if
$$\dim {\rm Hom}_Q(M,X)\le \dim {\rm Hom}_Q(P, X) + 1$$
for every non-injective indecomposable representation $X$.
\end{conjA} 

If Conjecture \ref{conj:hom-dim} holds, $M^2$ is described as a representation with maximal possible values of $\dim {\rm Hom}_Q(N,X)$ not exceeding $\dim {\rm Hom}_Q(P, X) + 1$.

Finally, we consider the Pl\"ucker algebras (the homogeneous coordinate rings)
${\rm Pl}(M)$, $M\in {\rm Rep}_\bd$ of the quiver Grassmannians $\Gr_{\dim P}(M)$.
These algebras are defined as the quotients of the polynomial ring in all Pl\"ucker
variables by the defining ideal  ${\mathcal I}(M)$ responsible for the reduced
scheme structure of our quiver Grassmannians. Let  ${\rm Pl}_\bfm(M)$ be the homogeneous 
components of ${\rm Pl}(M)$. 

\begin{conjA}\label{conj:Pl}
Let $Q$ be any Dynkin quiver with any orientation and let $P$ and $I$ be 
multiplicity free (i.e. each indecomposable summand shows up at most once).
Then for $M\in{\rm Rep}_\bd$ one has 
\[
\dim {\rm Pl}_\bfm(M) = \dim {\rm Pl}_\bfm(M^0) \text{ for all } \bfm
\]
if and only if $\Gr_{\dim P}(M)$ is of dimension $\langle \dim P,\dim I\rangle$.  
For arbitrary $M$ one has $\dim {\rm Pl}_\bfm(M) \ge \dim {\rm Pl}_\bfm(M^0)$
for all $\bfm$.
\end{conjA}

We note that Conjecture \ref{conj:Pl} does not hold for arbitrary projective $P$ and injective $I$ (even in type $A$).
At the same time, the conjectural flatness of the equi-dimensional (of minimal dimension) family of quiver Grassmannians 
$\Gr_{\dim P}(M)$ implies that the Euler characteristic of the natural line bundles
induced by the natural maps to $\bP(\Lambda^{(\dim P)_i}M_i)$ are constant along the
fibers. It is tempting to conjecture that for the multiplicity free $P$ and $I$ 
the components of the Pl\"ucker algebra are isomorphic to the zero cohomology groups 
of the corresponding line bundles and that the higher cohomology vanish. However,
at the moment we do not have enough evidence to conjecture that this indeed holds  true.  

Finally. let us outline couple of possible further directions. First, the principal quiver Grassmannians for equioriented type $A$ quivers enjoy many nice topological and combinatorial properties \cite{Bi14,CI20,CFFFR20,CFR13-1,Fe12-2}. 
It would be interesting to extend the study to the case of general quivers.
Second, a quiver Grassmannian of subrepresentations of $M$ is naturally acted upon by the group of automorphisms of $M$.
The description of this action and of the induced action on the homogeneous coordinate rings will lead to a better understanding of
geometric and algebraic properties of the quiver Grassmannians
(see \cite{Ar11,CFR17,Fe23,FFL11,HT99} for some partial results).

Our paper is organized as follows. In section \ref{sec:setup} we fix the notation and recall the basic definitions from the theory of quivers. 
In section \ref{sec:rss} we formulate the conjectural description of the reduced scheme structure of the principal quiver Grassmannians.
In section \ref{sec:irr} we formulate a conjecture on the locus of the irreducible principal quiver Grassmannians of the minimal possible dimension. 
Sections \ref{sec:flat} and \ref{sec:hom} treat the case of (possibly reducible) 
principal quiver Grassmannians of the minimal possible dimension.
In Section \ref{sec:hcr} we discuss the conjectural properties of the homogeneous coordinate rings.
Section \ref{sec:comp} contains a brief description of the software used to compute examples supporting conjectures given 
in the paper.

\section*{Acknowledgments}
We are grateful to Markus Reineke for useful discussions.
EF was partially supported by the ISF grant 493/24.

\section{The setup}\label{sec:setup}
In this section we recall main definitions from the theory of quiver representations and quiver Grassmannians \cite{ASS06,CB92,CR00,Schi14}. 
\subsection{Generalities}
We fix an algebraically closed field ${\mathbf k}$ of characteristic zero.
Let $Q$ be a Dynkin quiver with the set of vertices $Q_0$ and the set of arrows $Q_1$. 
For an arrow $\alpha$ we write $s(\alpha)$ for its source and $t(\alpha)$ for its target. A representation of $Q$ is a collection of vector spaces $M_i$, $i\in Q_0$ 
and a collection of linear maps $M_\alpha: M_{s(\alpha)}\to M_{t(\alpha)}$.

The simple modules of $Q$ are labeled by the vertices of $Q$. For $i\in Q_0$
we denote by $S_i$ the simple $1$-dimensional module supported at the vertex $i$.
There are finitely many iso-classes of indecomposable representations of $Q$;
these iso-classes are in bijection with positive roots of the simple Lie algebra 
whose Dynkin diagram is $Q$. In what follows for a vertex $i$ we denote by
$P_i$ the corresponding indecomposable projective module and by $I_i$ the
corresponding  indecomposable injective  module. 
In particular, one has  
the canonical embedding  $S_i\to I_i$ and the canonical projection $P_i\to S_i$.  
We note that $(P_i)_j$, $i,j\in Q_0$
is non-zero if and only if there exists a path from $i$ to $j$ in $Q$ and 
$(I_i)_j$, $i,j\in Q_0$
is non-zero if and only if there exists a path from $j$ to $i$.

For a dimension vector $\bd\in\bZ_{\ge 0}^{Q_0}$ let ${\rm Rep}_\bd$ be the representation 
space of $Q$ of dimension $\bd$:
\[
{\rm Rep}_\bd = \bigoplus_{\alpha\in Q_1} {\rm Hom} (\bk^{d_{s(\alpha)}},\bk^{d_{t(\alpha)}}).
\] 
The space ${\rm Rep}_\bd$ is naturally acted upon (via the base change) by the group $G_\bd$, 
which is the product of the general linear groups $GL_{d_i}$ for $i\in Q_0$.
The orbits of this action parametrize the isomorphism classes of the $\bd$-dimensional 
representations of $Q$. We denote by ${\rm Rep}(Q)$ the set of all representations 
of the quiver $Q$.

For a point $M\in {\rm Rep}_\bd$ let $\eO(M)\subset {\rm Rep}_\bd$ be the orbit $G_\bd M$.
The dimension and codimension of an orbit are computed by the following formulas:
\[
\dim \eO(M) = \dim G_d - \dim{\rm End}_Q(M),\  
{\rm codim}_{{\rm Rep}_\bd} \eO(M) = \dim {\rm Ext}^1_Q(M,M).
\]
In particular, an orbit $\eO(M)$ is dense in the representation space if and only if $M$
is rigid, i.e. ${\rm Ext}^1_Q(M,M)=0$. For a Dynkin quiver $Q$ there is only finitely many
isoclasses of $Q$ modules of a given dimension. Hence, there exists a unique open dense 
orbit. 

We say that $M$ degenerates to $N$ 
if $\overline{\eO(M)}\supset \eO(N)$; we write $M\le N$. One has the following theorem of
Bongartz \cite{Bo96}: $M\le N$ if and only if 
\[
\dim {\rm Hom}_Q(M,X)\le \dim {\rm Hom}_Q(N,X)\quad  \forall X\in {\rm Rep}(Q)
\]
or, equivalently, if 
\[
\dim {\rm Hom}_Q(X,M)\le \dim {\rm Hom}_Q(X,N)\quad \forall X\in {\rm Rep}(Q).
\]
The degeneration order induces the structure of a poset on the set of isoclasses of
$Q$-modules of a fixed dimension $\bd$: we say that the isoclass of $M$ is less than or 
equal to the isoclass of $N$ if $M$ degenerated to $N$. We denote this poset by $\Gamma_\bd$.
In particular, $\Gamma_\bd$ has a unique minimal element corresponding to the open dense
orbit.

\begin{figure}[htbp]
	\centering
	\includegraphics[width=0.3\linewidth]{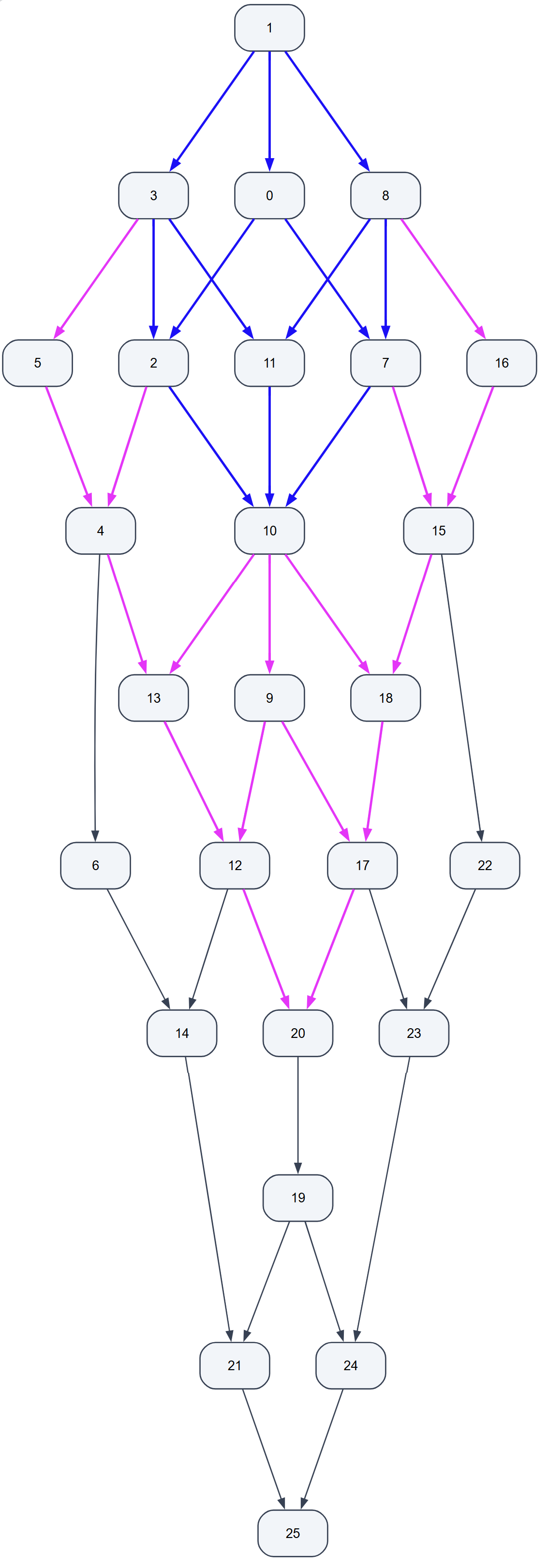}
	\caption{Degeneration picture for $Q=\bullet \rightarrow \bullet\leftarrow \bullet$,
		$\bd = (3,4,3)$.
	}
\end{figure}

\begin{rem}
Let $Q$ be the equioriented type $A$ quiver. Then the Bongartz criterion can be
explicitly formulated as follows:  $M$ degenerates to $N$ if and only if, for any path in $Q$, the rank
of the map in $M$ corresponding this path is no smaller than the 
corresponding rank of the map in $N$ \cite{AdF84}.
The general type $A$ case is worked out on \cite{AdF85}. 
\end{rem}	

For two vectors ${\bf a},{\bf b}\in{\mathbb Z}^{Q_0}$ we denote by $\langle {\bf a}, {\bf b}\rangle$
the value of the Euler form on these vectors given by
\[
\langle {\bf a}, {\bf b}\rangle = \sum_{i\in Q_0} a_ib_i - \sum_{\alpha\in Q_1} a_{s(\alpha)}b_{t(\alpha)}.
\]  
For representations $M\in {\rm Rep}_{\bf a}$ and $N\in {\rm Rep}_{\bf b}$ one has
\[
\langle {\bf a}, {\bf b}\rangle = \dim \Hom_Q(M,N) - \dim {\rm Ext}^1_Q(M,N).
\]

\subsection{Quiver Grassmannians}
Let $M$ be a $\bd$-dimensional representation of $Q$ and let 
${\bf e}\in\bZ_{\ge 0}^{Q_0}$
be a dimension vector such that $e_i\le d_i$ for all $i\in Q_0$. The universal quiver Grassmannian
$\Gr_{\bf e}(\bd)$ sits inside the product ${\rm Rep}_\bd\times\prod_{i\in Q_0} \Gr_{e_i}(\bk^{d_i})$
and consists of collections $\{f_\alpha\}_{\alpha\in Q_1}, (U_i)_{i\in Q_0}$ such that 
$f_\alpha(U_{s(\alpha)})\subset U_{t(\alpha)}$ for all $\alpha\in Q_1$. One has two natural
projections    
\[
{\rm Rep}_\bd \longleftarrow  \Gr_{\bf e}(\bd) \longrightarrow   \prod_{i\in Q_0} \Gr_{e_i}(\bk^{d_i});
\]
for $M\in {\rm Rep}_\bd$ we denote by $\Gr_{\bf e}(M)$ the fiber of the projection $\Gr_{\bf e}\to {\rm Rep}_\bd$ over the point $M$.

Let us fix a projective representation $P$ and an injective representation $I$;
thus $P$ is a direct sum of the modules $P_i$ with certain multiplicities and 
$I$ is a direct sum of the modules $I_i$ with certain multiplicities. Let 
$$\bd=\dim P + \dim I.$$ 
We are interested in quiver Grassmannians $\Gr_{\dim P}(M)$ for $M\in {\rm Rep}_\bd$
being $\bd$-dimensional
$Q$-module. In other words, we want to study the fibers of the universal quiver Grassmannian
$\Gr_{\dim P}(P + I)$ (over the representation space ${\rm Rep}_\bd$). Following 
suggestion of Markus Reineke we call these quiver Grassmannians PrIncipal ({\it Pr} for  projectives and {\it In} for injectives).

Since $Q$ is Dynkin, the group $G_\bd$ has a finite number of orbits in ${\rm Rep}_\bd$.
Hence there exists a unique open orbit $\eO(M^0)$ (for a rigid representation $M^0$), which
degenerates to any other orbit and a special fiber $\Gr_{\dim P}(M^0)$ of $\Gr_{\dim P}(\bd)$.  
Yet another special fiber is given by $\Gr_{\dim P}(P\oplus I)$ (since
$\bd=\dim P + \dim I$).  

\begin{rem}
For any $M\in {\rm Rep}_\bd$ the quiver Grassmannian $\Gr_{\dim P}(M)$ is non-empty. 
In fact, one knows \cite{Scho92} that this statement is implied by the inequalities
$\langle \ba, \dim I\rangle\ge 0$ for any dimension vector $\ba$, which holds true because $I$ is injective. Furthermore, if $Q$ is of type $A$, then
it follows from \cite{CEFR21,CEFF22}
that for any $M\in {\rm Rep}_\bd$ the quiver Grassmannian
$\Gr_{\dim P}(M)$ is connected. 
More precisely, one shows that all quiver Grassmannians for the type $A$ quivers are connected. It is tempting to conjecture that
the same holds true for other Dynkin quivers as well: all the type $D$ quiver Grassmannians  $\Gr_{\dim P}(M)$ we've seen in our computations are indeed connected.  
\end{rem}

The following statements were proved in \cite{CFR12}.
\begin{itemize}
\item $\Gr_{\dim P}(M^0)$ and $\Gr_{\dim P}(P\oplus Q)$ are irreducible, 
\item $\dim \Gr_{\dim P}(M^0) = \dim \Gr_{\dim P}(M^0) = \langle \dim P,\dim I\rangle$,
\item  $\Gr_{\dim P}(P\oplus Q)$ is a flat degeneration of $\Gr_{\dim P}(M^0)$,
\item if $M$ degenerates to $P\oplus Q$, then $ \Gr_{\dim P}(M)$ is irreducible of dimension 
$\langle \dim P,\dim I\rangle$.
\end{itemize}

\subsection{Equioriented type  A case}
Let $Q$ be the equioriented quiver of type $A_n$.  We denote the vertices of $Q$ by $1,2,\dots,n$ with arrows 
given by $i\to i+1$.  The indecomposable representations are labeled by pairs $i,j$ with
$1\le i\le j\le n$; we denote the corresponding representation by $U(i,j)$.
In particular, all  non-trivial components of $U(i,j)$ are one-dimensional and the 
support of $U(i,j)$ consists of vertices between $i$ and $j$. One has $S_i=U(i,i)$,
$P_i = U(i,n)$, $I_i = U(1,i)$. The representation $M^0$ is equal to $U(1,n)^{\oplus n+1}$. 

Let us consider the special case $P=\bigoplus_{i=1}^n P_i$, $I=\bigoplus_{i=1}^n I_i$.
Then $P$ is isomorphic (as a $Q$-module) to the path algebra and $I$ to its dual.
One has (recall $\bd = \dim P +\dim I$)
\[
\dim A = (1,\dots,n),\ \dim A^* = (n,\dots,1), \ \bd = (n+1,\dots,n+1).
\]
The scalar product $\langle\dim A, \dim A^*\rangle$ is equal to $n(n+1)/2$. One observes 
that $\Gr_{\dim A}(M^0)$ is isomorphic to the classical full flag variety for the group
$SL_{n+1}$ (see e.g. \cite{Fu97}) and  $\Gr_{\dim A}(P\oplus Q)$ is isomorphic to the PBW degenerate 
flag variety for the group $SL_{n+1}$ (see e.g. \cite{Fe11,Fe12-1,Fe23}).

We also introduce one more representation of dimension $\bd$: 
\[
M^2 = \bigoplus_{i=1}^n U(i,n) \oplus \bigoplus_{i=1}^n U(i,i) \oplus \bigoplus_{i=2}^n U(1, i-1).
\]
The following  is proved in \cite{CFFFR17} (see also \cite{CFFFR20}).
\begin{itemize}
\item $\Gr_{\dim A}(M)$ is irreducible of dimension $n(n+1)/2$ if and only if $M$ degenerates to $M^1$,
\item $\Gr_{\dim A}(M)$ is of dimension $n(n+1)/2$ if and only if $M$ degenerates to $M^2$.
\end{itemize}

\begin{rem}
The results of \cite{CFFFR17} have to do with flat and flat irreducible loci of the universal quiver Grassmannian $\Gr_{\dim A}(\bd)$.
\end{rem}

\section{Reduced scheme structure}\label{sec:rss}
For a number $n\in\bZ_{>0}$ we write $[n]$ for the set $\{1,\dots,n\}$.

Let us denote the components of the dimension vector $\dim P$ by $(a_i)_{i\in Q_0}$.
A quiver Grassmannian $\Gr_{\dim P}(M)$ sits inside the product of classical Grassmannians
$\Gr_{a_i}(M_i)$, $i\in Q_0$. Using the standard Pl\"ucker embeddings of the Grassmann
varieties one arrives at the embedding
\[
\Gr_{\dim P}(M) \subset \prod_{i\in Q_0} \Gr_{a_i}(M_i) \subset \prod_{i\in Q_0} \bP(\Lambda^{a_i} (M_i)).
\]   
We call the composition of these embedding the Pl\"ucker embedding for quiver Grassmannians. 
Let $\Delta^{(i)}_J$, $J\in \binom{[d_i]}{a_i}$, be the Pl\"ucker coordinates in the wedge space
$\Lambda^{a_i} (M_i)$. We denote by the same symbols the homogeneous coordinates on the 
projective space $\bP(\Lambda^{a_i} (M_i))$. Let ${\mathcal I}(M)$ be the multi-homogeneous ideal in the 
polynomial ring in all the variables $\Delta^{(i)}_J$ (for all $i\in Q_0$)
consisting of all multi-homogeneous polynomials vanishing on $\Gr_{\dim P}(M)$. In order to give a (conjectural)
explicit description of ${\mathcal I}(M)$ we introduce some notation.

Let us fix a basis $v^{(i)}_p$, $1\le p\le d_i$ in each component $M_i$ of the representation $M$.   
For an arrow $\alpha: i\to j$, $\alpha\in Q_1$ let $m_{\alpha,p,q}$ be the matrix of the map $M_\alpha$
written in the fixed bases. For a number $p\in [d_i]$ and a subset $I\subset [d_i]$ we write 
$\epsilon(p,I)=\#\{p'\in I:\ p'\le p\}$. The following relations are labeled by an arrow 
$\alpha:i\to j$, a set $I\subset [d_i]$ of cardinality $a_i-1$ and a set $J\subset [d_j]$ of cardinality $[d_j+1]$:
\begin{equation}\label{eq:relationsarrows}
R(\alpha,I,J) = \sum_{\substack{p\in [d_i]\setminus I\\ q\in J}} (-1)^{\epsilon(p,I)+\epsilon(q,J)}
m_{\alpha,p,q}\Delta^{(i)}_{I\cup p}\Delta^{(j)}_{J\setminus \{q\}}.
\end{equation}
It is shown in \cite{LW19} that the relations above cut out $\Gr_{\dim P}(M)$ pointwise inside the product
of projectivized wedge powers. Similarly to relations  \eqref{eq:relationsarrows} one gets a relation 
$R(\pi,I,J)$
for any path $\pi$ in $Q$ with $s(\pi)=i$, $t(\pi)=j$ and  $I\in\binom{[d_i]}{a_i-1}$,
$J\in\binom{[d_j]}{a_j+1}$:
\begin{equation}\label{eq:relationspaths}
	R(\pi,I,J) = \sum_{\substack{p\in [d_i]\setminus I\\ q\in J}} (-1)^{\epsilon(p,I)+\epsilon(q,J)}
	m_{\pi,p,q}\Delta^{(i)}_{I\cup p}\Delta^{(j)}_{J\setminus \{q\}},
\end{equation}
with $m_{\pi,p,q}$ being matrix coefficients of the map $M_\pi$ (the composition of maps corresponding
to the arrows in $\pi$). 

Recall the ideal ${\mathcal I}(M)$ defining the reduced scheme structure of 
$\Gr_{\dim P}(M)$.
We put forward the following conjecture:
\begin{conj}\label{conj:relations}
For any $M\in {\rm Rep}_\bd$ 
the ideal ${\mathcal I}(M)$ is generated by the quadratic relations $R(\pi,I,J)$ for all path $\pi$ in $Q$ and
all  $I\subset [d_{s(\pi)}]$, $\# I=a_{s(\pi)-1}$ and $J\subset [d_{t(\pi)}]$, $\# J= a_{t(\pi)}+1$. 	
\end{conj}	     

In particular, the ideal generated by all the relations  $R(\pi,I,J)$ is saturated (with respect to 
each group of Pl\"ucker variables $\Delta^{(i)}_\bullet$) and is prime.   
We note that Conjecture \ref{conj:relations} is far from being true for general quiver
Grassmannians.

\section{Minimal dimension locus: irreducible Grassmannians}
\label{sec:irr}

Recall that $\dim \Gr_{\dim P}(M^0) = \dim \Gr_{\dim P}(P\oplus Q) = \langle \dim P, \dim I\rangle$ and
both varieties are irreducible. Moreover, if a representation $M\in{\rm Rep}_\bd$ degenerates 
to $P\oplus Q$, then the quiver Grassmannian $\Gr_{\dim P}(M)$ is also irreducible and of the 
same dimension. 

Let's denote by $\Gamma_\bd(1)$ the subset of the degeneration poset $\Gamma_\bd$ that consists of representations $M\in {\rm Rep}_\bd$ whose Grassmanian $\Gr_{\dim P}(M)$ is irreducible of dimension $\langle \dim P, \dim I\rangle$. 

Clearly, the lower ideal of $\Gamma_\bd$ consisting of representations $M$ that degenerate to $P\oplus Q$ is contained in $\Gamma_\bd(1)$. It turns out that they often coincide.    

\begin{example}
Let $Q$ be the $A_3$ quiver of the form $\bullet \rightarrow \bullet \leftarrow \bullet$. 
Let $P$ and $I$ be the direct sums of all projective (resp., injective) indecomposable 
representations. Then the poset $\Gamma_\bd(1)$ is visualized at Figure  \ref{fig:A3sink1} below. Its single sink is exactly $P\oplus Q$.
\end{example}
\begin{figure}[htbp]
	\centering
	\includegraphics[width=0.3\linewidth]{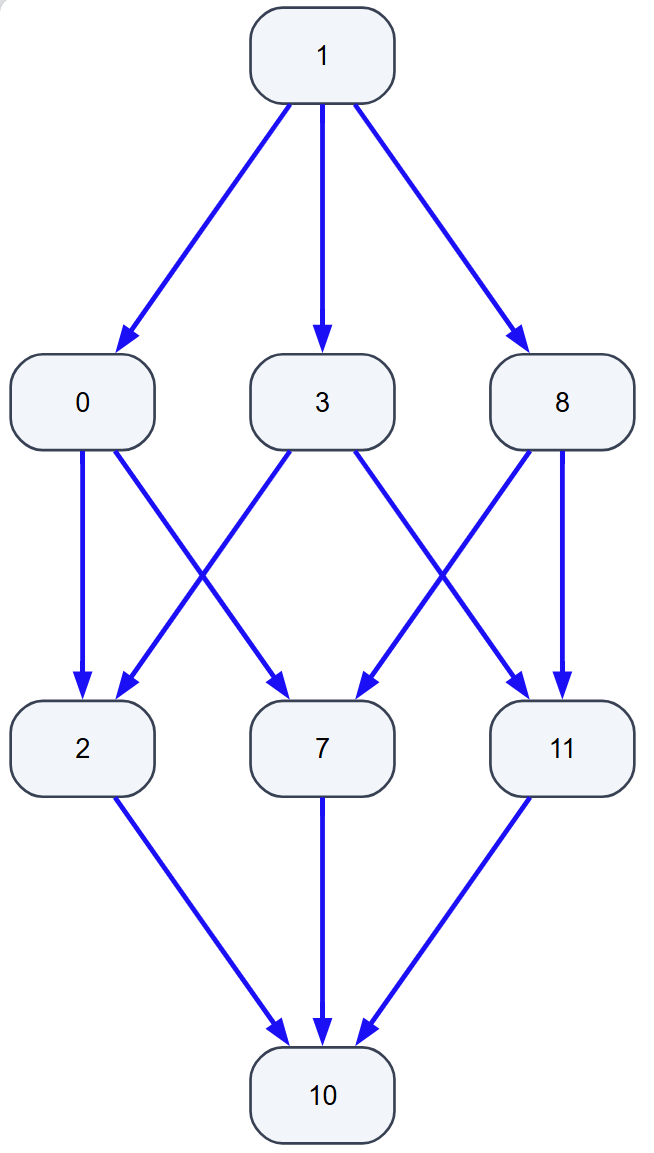}
	\caption{$\Gamma_\bd(1)$ for $Q=\bullet \rightarrow \bullet\leftarrow \bullet$,
		$P=P_1\oplus P_2\oplus P_3$, $I=I_1\oplus I_2\oplus I_3$.
	}
	\label{fig:A3sink1}
\end{figure}

\begin{conj}
There exists a representation $M^1$ such that for a representation $M\in{\rm Rep}_\bd$ the quiver Grassmannian $\Gr_{\dim P}(M)$ is irreducible of
dimension $\langle \dim P, \dim I\rangle$ if and only if $M$ degenerates to $M^1$ (i.e.
$M\in\Gamma_\bd(1)$ if and only if $M^0\le M\le M^1$). 
If, furthermore, no component of either $\dim P$ and $\dim I$ is zero, $M^1 = P\oplus I$.
\end{conj} 

The statement is proved in \cite{CFFFR17,CFFFR20} for certain $P$ and $I$ for the equioriented type $A$ quivers.

The condition that $\dim P$ and $\dim I$ do not have zero components seems to be crucial. In all the cases we checked, where some of the components vanish
(both for $A_n$ and $D_n$ quivers), $M^1$ was different from $P\oplus I$ and thus $\{M:\ M^0\le M\le P\oplus I\}$ was a strict subgraph of $\Gamma_\bd(1)$.

\begin{example}
Let $Q$ be the $A_3$ quiver of the form $\bullet \rightarrow \bullet \leftarrow \bullet$. 
Let $P = P_1\oplus P_2\oplus P_3$ and $I = I_1\oplus I_3$. Then $\bd = (2, 3, 2)$, $\dim I = (1, 0, 1)$, and $\dim P = (1, 3, 1)$. 

In this case, whatever the representation $M$ is, any one-dimensional subspace of $M_1$ and any one-dimensional subspace of $M_3$ define a subrepresentation of dimension $\dim P$. Thus, $\Gr_{\dim P}(M)\cong\mathbb{P}^1\times\mathbb{P}^1$ for every $M$. So, $\Gamma_\bd(1) = \Gamma_\bd$ and $M^1 = 2S_1\oplus 3S_2\oplus 2S_3$ which is not $P\oplus I$.
\end{example}

\begin{example}
Let $Q$ be the $A_4$ quiver of the form $\bullet \rightarrow \bullet \leftarrow \bullet \leftarrow \bullet$. 
Let $P = \oplus_{i=1}^4 P_i$ and $I = I_1\oplus I_3\oplus I_4$. Then $\bd = (2, 4, 3, 3)$, $\dim I = (1, 0, 1, 2)$, and $\dim P = (1, 4, 2, 1)$.

In this case $\Gamma_\bd(1) \ne \Gamma_\bd$, but
$$M^1 = 2U(1,1) + 4U(2,2) + U(3, 3) + 2U(3, 4) + U(4, 4),$$
which differs from
\[
P\oplus I = U(1, 1) + U(1, 2) + U(2, 2) + U(2, 3) + U(2, 4) + U(3, 4) + U(4, 4).
\]
\end{example}

In the degenerate examples that we studied, the representation $M^1$ can be constructed from $P\oplus I$ in the following way. Let's call a vertex $k$ \textit{deficient} if $(\dim P)_k = 0$ or $(\dim I)_k = 0$. For every direct summand of $P\oplus I$, ``split'' it at every deficient vertex by vanishing the maps incidental to this vertex. The resulting sum will be $M^1$.

In the previous $A_4$ example, the deficient vertex is $2$, and we have four summands $M$ in $P\oplus I$ such that $M_2\ne 0$. Let's split those supported not only at $2$:
\begin{align*}
U(1,2)&\mapsto U(1,1) + U(2, 2),\\
U(2,3)&\mapsto U(2,2) + U(3, 3),\\
U(2,4)&\mapsto U(2,2) + U(3, 4).
\end{align*}
The resulting sum
\begin{align*}
&U(1, 1) + (U(1,1) + U(2, 2)) + U(2, 2) + (U(2,2) + U(3, 3)) +\\ 
& (U(2,2) + U(3, 4)) + U(3, 4) + U(4, 4)
\end{align*}
is exactly $M^1$.

\section{Minimal dimension locus}\label{sec:flat}
In this section we study the locus of representations $M$ such that the quiver Grassmannian 
$\Gr_{\dim P}(M)$ is of the expected (minimal possible) dimension
$\langle\dim P,\dim I\rangle$, but may have more than one irreducible component.
We denote by $\Gamma_\bd(2)$ the poset of isoclasses of such representations $M$. 
Clearly, if $N\in  \Gamma_\bd(2)$ and $M\le N$, then $M\in  \Gamma_\bd(2)$ and hence
$\Gamma_\bd(2)\subset \Gamma_\bd$ is a lower ideal, which contains a subideal $\Gamma_\bd(1)$.
 
\begin{example}
Let $Q$ be the $A_3$ quiver of the form $\bullet \rightarrow \bullet \leftarrow \bullet$. 
Let $P$ and $I$ be the direct sums of all projective (resp., injective) indecomposable 
representations. Then the poset $\Gamma_\bd(2)$ is visualized at Figure 
\ref{fig:A3_sink_2} below.
\end{example}

\begin{figure}[htbp]
	\centering
	\includegraphics[width=0.3\linewidth]{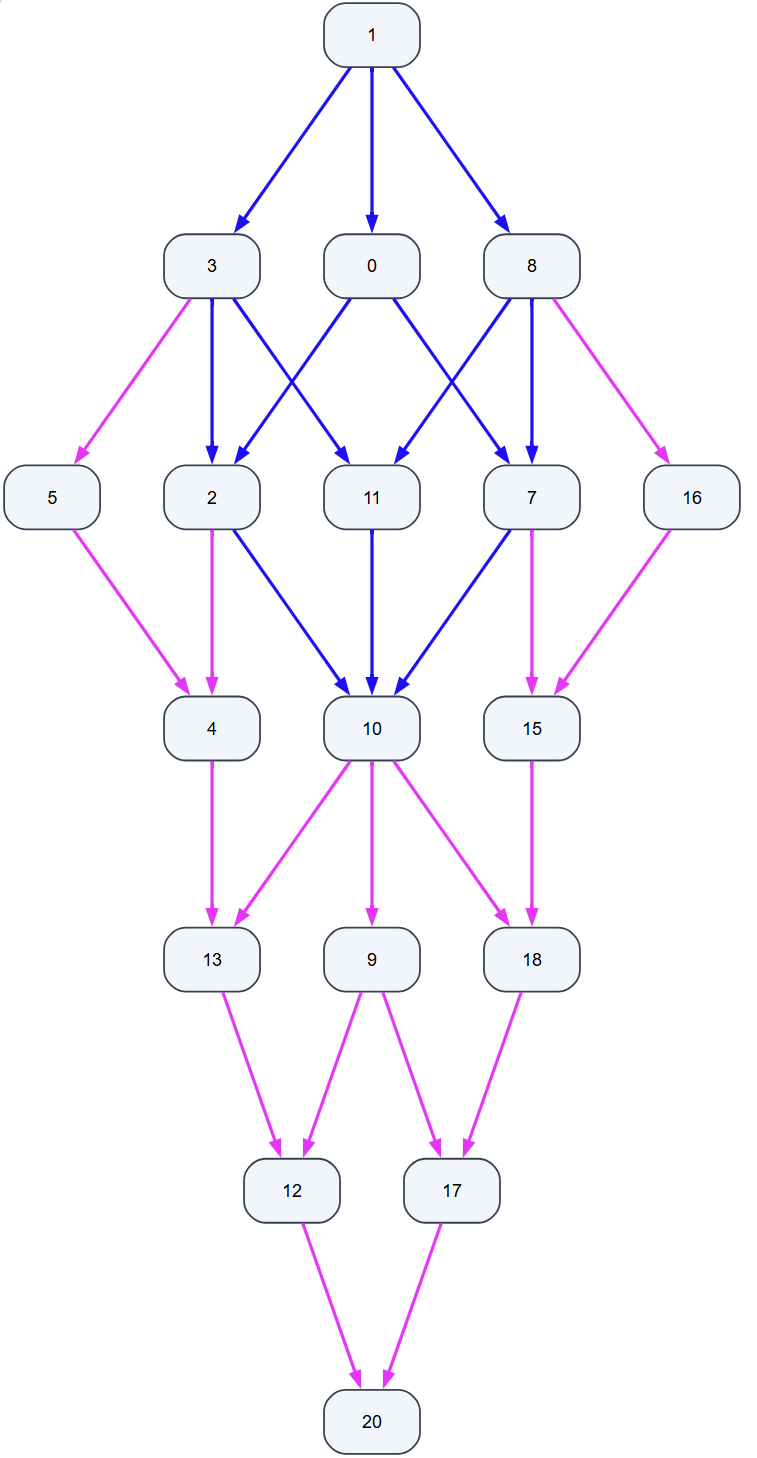}
	\caption{$\Gamma_\bd(2)$ for $Q=\bullet \rightarrow \bullet\leftarrow \bullet$,
		$P=P_1\oplus P_2\oplus P_3$, $I=I_1\oplus I_2\oplus I_3$.
	}
	\label{fig:A3_sink_2}
\end{figure}

\subsection{Type A}
Let $Q$ be a type $A$ quiver with arbitrary orientation. We put forward the following conjecture.

\begin{conj}
There exists a unique representation
$M^2\in {\rm Rep}_\bd$ such that $\Gr_{\dim P}(M)$ is of dimension $\langle \dim P, \dim I\rangle$ 
if and only if $M$ degenerates to $M^2$. If $M$ degenerates to $M^2$, then $\Gr_{\dim P}(M)$
is equidimensional. 
\end{conj} 

The conjecture above can be rephrased in the following way: the poset $\Gamma_\bd(2)$ consists 
of elements $M$ such that $M^0\le M\le M^2$ for certain element $M^2\in\Gamma_\bd$. 

In general we expect that the following conjecture holds true.
\begin{conj}
If $\dim \Gr_{\dim P}(M)=\langle \dim P, \dim I\rangle$, then $\dim \Gr_{\dim P}(M)$ 
is equidimensional.
\end{conj} 

Recall (see \cite{CFFFR17}) that for the equioriented type $A_n$ quiver and $P=\bigoplus_{i\in Q_0} P_i$,
$I=\bigoplus_{i\in Q_0} I_i$
the representation $M^2$ as above does exist 
and is isomorphic to
\begin{equation}\label{eq:M2A}
M^2 = P\oplus S\oplus I/S, \quad S=\bigoplus_{i\in Q_0} S_i.
\end{equation}
The corresponding 
quiver Grassmannian $\Gr_{\dim P}(M^2)$ is equidimensional, of expected dimension $n(n+1)/2$ and has 
the Catalan number irreducible components.

Let us adjust the description \eqref{eq:M2A} above to the conjectural answer in type $A$ for general orientation. We decompose 
$Q$ into the union of several equioriented type $A$ quivers
$Q=Q(1)\cup\dots\cup Q(r)$ such that for each $\ell$ the quivers
$Q(\ell)$ and $Q(\ell+1)$ have a unique common vertex $i(\ell)\in Q_0$. Let $M^2(\ell)$ be the deepest representation 
\eqref{eq:M2A} corresponding to the equioriented quiver $Q(\ell)$. In particular, each $M^2(\ell)$ contains summands of the form
$2S(i)$ for $i$ being the leftmost and the rightmost vertices of $Q(\ell)$. 

We put forward the following conjectures: 

\begin{conj}\label{conj:M^2-split}
Let $Q$ be of type $A$ with the decomposition $Q=\cup_{\ell=1}^r Q(\ell)$, $Q(\ell)\cap Q(\ell+1)=i_\ell$. Let also $P = \bigoplus_{j\in Q_0} P_j$ and $I = \bigoplus_{j\in Q_0} I_j$. Then 
\[
M^2 = \bigoplus_{\ell=1}^r M^2(\ell) - 2 \bigoplus_{a=1}^{r-1} S_{i_a},
\]
where the subtraction means that we remove two summands of the form $S_{i_a}$, $1\le a <r$ from the decomposition of  
$\bigoplus_{\ell=1}^r M^2(\ell)$ into the direct sum of indecomposable representations.
\end{conj}

\begin{conj}\label{conj:M^2-higher}
Let $P = \bigoplus_ju_j^PP_j$, $I = \bigoplus_ju_j^II_j$ with all $u_j^P, u_j^I\ne 0$. Let also $M^2(P,I)$ be the corresponding deepest representation with minimal dimension Grassmannian. Then
$$M^2(P\oplus P', I\oplus I) = M^2(P, I)\oplus P'\oplus I'.$$
\end{conj}

In the case if some $u_j^P$ or $u_j^I$ are zero, the exact form of $M^2$ is not entirely clear for us yet. In some cases Conjecture \ref{conj:M^2-higher} still holds, in some it doesn't.

We provide some concrete examples below.

\subsection{Zig-zag case}
We continue using the notation $U(i,j)$ for the indecomposable $A_n$ module
supported on vertices between $i$ and $j$ (for  arbitrary orientation).
In the following Conjecture we consider quiver $A_n$ with alternating orientation. 

\begin{conj}
Let $Q$ be of type $A_n$ with alternating orientation (a ``zig-zag quiver''). Then $M^2$ is written as 
$\bigoplus_{i=1}^{n-1} U(i,i+1)\oplus 2\bigoplus_{i\in Q_0} S_i$. The quiver Grassmannian
$\Gr_{\dim P}(M^2)$ has $2^{n-1}$ irreducible components each of dimension $2n-1$. 
\end{conj} 

\subsection{$A_5$ example}
Let us consider an $n=5$ example. Let $Q$ be as follows:
$1 \rightarrow 2 \leftarrow 3 \leftarrow 4 \leftarrow 5$. Then the representation $M^2$ 
is given by the following formula:
\begin{multline*}
M^2 = 2U(1,1) + U(1,2) + 2U(2,2) +  U(2,3) + U(3,3) + U(2,4) +\\ U(4,4) + 
U(2,5) + U(3,5) + U(4,5) + 2U(5,5).
\end{multline*}
We note that the right hand side can be rewritten as    
\begin{multline*}
\Bigl(2U(1,1) + U(1,2) + 2U(2,2)\Bigr)  + \Bigl( 2U(2,2) + U(2,3) + U(3,3) + U(2,4) +\\ U(4,4) + 
	U(2,5) + U(3,5) + U(4,5) + 2U(5,5)\Bigr) - 2U(2,2),
\end{multline*}
where the sums in brackets correspond to the representations $M^2$ for the quivers 
$1 \rightarrow 2$ and $2 \leftarrow 3 \leftarrow 4 \leftarrow 5$ (with $- 2U(2,2)$ coming as gluing at vertex $2$). 

\subsection{Type D}
Let $Q$ be the following $D_4$ quiver

\begin{picture}(100,110)     
	\put(92,50){\texttt{1}}    
	\put(100,54){\vector(1,0){30}} 
	\put(134,50){\texttt{2}}    
	\put(142,58){\vector(1,1){30}} 
	\put(142,46){\vector(1,-1){30}} 
	\put(174,90){\texttt{3}}    
	\put(174,10){\texttt{4}}    
\end{picture}

\noindent Computer experiments support the following proposition.

\begin{prop}
Let $Q$ the a quiver with vertices $1,2,3,4$ and arrows $1\to 2$, $2\to 3$ and $2\to 4$. 	
Let $P=\bigoplus_{i\in Q_0} P_i$, $I=\bigoplus_{i\in Q_0} I_i$. 
Then there exist three representations $M^2(1)$, $M^2(2)$, $M^2(3)$  in ${\rm Rep}_\bd$ such that
$\dim \Gr_{\dim P}(M) =\langle \dim P, \dim I\rangle$ if and only if $M$ degenerates to one of 
$M^2(i)$, $i=1,2,3$. None of the representations $M^2(1)$, $M^2(2)$, $M^2(3)$ degenerate one to another.	
\end{prop} 

Let us provide more details. One has the following explicit formulas:
\[
\bd = (5,5,4,4),\ \dim P = (1,2,3,3),\ \langle \dim P, \dim I \rangle = 9.  
\]
In order to write down the representations $M^2(i)$ explicitly, we denote by $V(x_1,x_2,x_3,x_4)$
the indecomposable $Q$ module of dimension $(x_1,x_2,x_3,x_4)$. Then one has 
\begin{multline*}
M^2(1) = 2V(1,0,0,0) \oplus V(1,1,0,0) \oplus V(0,1,0,0) \oplus V(1,1,1,0) \oplus\\ 2V(0,0,1,0) \oplus V(0,1,0,1) \oplus 
2V(0,0,0,1) \oplus V(1,1,1,1),
\end{multline*}
\begin{multline*}
M^2(2) = 2V(1,0,0,0) \oplus V(1,1,0,0) \oplus V(0,1,0,0) \oplus V(0,1,1,0) \oplus\\ 2V(0,0,1,0) \oplus V(1,1,0,1) \oplus 
2V(0,0,0,1) \oplus V(1,1,1,1),
\end{multline*}
\begin{multline*}
M^2(3) = 2V(1,0,0,0) \oplus V(1,1,0,0) \oplus V(1,1,1,0) \oplus V(0,1,1,0) \oplus \\ 2V(0,0,1,0) \oplus V(1,1,0,1) \oplus 
V(0,1,0,1) \oplus 2V(0,0,0,1).
\end{multline*}
All three quiver Grassmannians $\Gr_{\dim P}(M^2(i))$  have $13$ irreducible components of dimension $9$.

Now let us consider the $D_4$ quiver with all vertices pointing to the central vertex

\begin{picture}(100,110)     
	\put(92,50){\texttt{1}}    
	\put(100,54){\vector(11,0){30}} 
	\put(134,50){\texttt{2}}    
	\put(172,88){\vector(-1,-1){30}} 
	\put(172,18){\vector(-1,1){30}} 
	\put(174,90){\texttt{3}}    
	\put(174,12){\texttt{4}}    
\end{picture}

\begin{prop}
Let $Q$ the a quiver with vertices $1,2,3,4$ and arrows $1\to 2$, $3\to 2$ and $4\to 2$. 	
Let $P=\bigoplus_{i\in Q_0} P_i$, $I=\bigoplus_{i\in Q_0} I_i$.
Then there exists a unique representations $M^2$ in ${\rm Rep}_\bd$ such that
$\dim \Gr_{\dim P}(M) =\langle \dim P, \dim I\rangle$ if and only if $M$ degenerates to $M^2$. 	
\end{prop} 

Let us provide the details in this case. 
One has 
\[
\bd = (3,5,3,3),\ \dim P = (1,4,1,1),\ \langle \dim P, \dim I \rangle = 7  
\]
and
\begin{multline*}
	M^2 = 2V(1,0,0,0) \oplus V(1,1,0,0) \oplus 2V(0,1,0,0) \oplus V(0,1,1,0) \oplus\\ 2V(0,0,1,0) \oplus V(0,1,0,1) \oplus 2V(0,0,0,1).
\end{multline*}
The quiver Grassmannian $\Gr_{\dim P}(M^2)$ has $8$ irreducible components of dimension $7$.

Finally, let us mention that for $Q$ of type $D_5$ of the form

\begin{picture}(100,110)     
\put(100,50){\texttt{1}}    
\put(108,54){\vector(1,0){30}} 
\put(142,50){\texttt{2}}    
\put(180,54){\vector(-1,0){30}} 
\put(184,50){\texttt{3}}    
\put(222,88){\vector(-1,-1){30}} 
\put(222,18){\vector(-1,1){30}} 
\put(224,90){\texttt{4}}    
\put(224,12){\texttt{5}}    
\end{picture}

\noindent one also gets three deepest representations $M^2$ such that the corresponding quiver Grassmannian
is of the expected (minimal possible) dimension.

\section{The homomorphism dimension criterion}\label{sec:hom}

The goal of this section is to formulate a criterion for a representation $N$ to belong to $\Gamma_\bd(2)$ in terms of dimensions of certain homomorphism spaces. 

We start with the following observation, which has to do with the description of $\Gamma_\bd(1)$:
\begin{itemize}
\item for every indecomposable injective $X = I_j$, $\dim {\rm Hom}_Q(M^1, X) = \bd_j$, 	\item for every indecomposable non-injective $X$, $\dim {\rm Hom}_Q(P\oplus I, X) = \dim {\rm Hom}_Q(P, X)$ 
\end{itemize}
(the first equality follows from the general fact that 
$\dim {\rm Hom}_Q(M, I_j) = \bd_j$ for every $M\in{\rm Rep}_\bd(Q)$, and the second holds true 
since there are no homomorphisms from $I$ to $X$).
Assume that the dimension vectors $\dim P$ and $\dim I$ never vanish.  In this case we expect (see section \ref{sec:irr}) that $M^1 = P\oplus Q$ and hence we can characterize 
$\Gamma_\bd(1)$ as the set of all representations $M$ such that, for every non-injective indecomposable representation $X$,
\[
\dim {\rm Hom}_Q(M,X)\le \dim {\rm Hom}_Q(P, X).
\]



It turns out that $M^2$ allows for a similar description.

\begin{conj}\label{conj:M^2-hom}
Let $P\oplus I$ contains every indecomposable projective and every indecomposable injective representation as a summand. Then a representation $M$ lies in $\Gamma_\bd(2)$ if and only if, for every non-injective indecomposable representation $X$,
\[
\dim {\rm Hom}_Q(M,X)\le \dim {\rm Hom}_Q(P, X) + 1.
\]
\end{conj} 

\begin{rem}
Dually, this can be formulated as: a representation $M$ lies in $\Gamma_\bd(2)$ if and only if, for every non-projective indecomposable representation $X$,
\[
\dim {\rm Hom}_Q(X,M)\le \dim {\rm Hom}_Q(X, I) + 1/
\]
\end{rem}

This provides a general way of describing $M^2$ for both $A$ and $D$ types. 
For $A$ type quivers, there exists a representation $M^2$ satisfying, for all indecomposable non-injective representations $X$,
\[
\dim {\rm Hom}_Q(M^2,X) = \dim {\rm Hom}_Q(P, X) + 1.
\]
For $D$ type quivers such a representation $M^2$ does not always exist. That's why in some cases we have three sinks $M^2(i)$ in $\Gamma_\bd(2)$. However, 
\[
\max_i\dim {\rm Hom}_Q(M^2(i),X) = \dim {\rm Hom}_Q(P, X) + 1
\]
for all indecomposable non-injective $X$ (assuming the conditions of Conjecture \ref{conj:M^2-hom} are satisfied).
Note that this perfectly aligns with Conjecture \ref{conj:M^2-higher}. Indeed, if the above equation holds, then
\[
\max_i\dim {\rm Hom}_Q(M^2(i)\oplus P'\oplus I',X) = \dim {\rm Hom}_Q(P\oplus P', X) + 1.
\]


If some of $P_j$ or $I_j$ are missing in $P\oplus I$, then Conjecture \ref{conj:M^2-hom} sometimes fails. 
For example, let $Q = D_4$ with the arrow orientation illustrated below, $P = P_1\oplus P_3\oplus P_4$, $I = I_1\oplus I_2\oplus I_3\oplus I_4$. 

\begin{picture}(100,110)     
	\put(92,50){\texttt{1}}    
	\put(100,54){\vector(1,0){30}} 
	\put(134,50){\texttt{2}}    
	\put(142,58){\vector(1,1){30}} 
	\put(142,46){\vector(1,-1){30}} 
	\put(174,90){\texttt{3}}    
	\put(174,10){\texttt{4}}    
\end{picture}

Then $\Gamma_\bd(2)$ has a single sink
\begin{align*}
M^2 = & 2V(1, 0, 0, 0) \oplus V(1, 1, 0, 0) \oplus V(0, 1, 0, 0) \oplus V(1, 1, 1, 0) \oplus\\ 
&2V(0,0,1,0) \oplus V(1, 1, 0, 1) \oplus 2V(0,0,0,1),
\end{align*}
such that
\[
\dim {\rm Hom}_Q(M^2, V(0,1,1,1)) = 4 = \dim {\rm Hom}_Q(P, V(0,1,1,1)) + 2.
\]
This case has one more peculiarity: the set of representations $N$ such that $\dim {\rm Hom}_Q(N, X) \le \dim {\rm Hom}_Q(P, X) + 1$, for all non-injective indecomposable $X$, has four sinks instead of one or three. At the same time, $M^1$ in this case coincides with $P\oplus I$.

On the other hand, if $Q$ and $I$ are the same  and $P = P_1\oplus P_2\oplus P_4$, the conjecture holds.

\section{Pl\"ucker algebras and line bundles}\label{sec:hcr}
Recall the the notation $a_i=\dim P_i$ and the
standard closed embedding $\Gr_{\dim P}(M)\hookrightarrow\prod_{i\in Q_0} \Gr_{a_i}(M_i)$. Each Grassmann
variety admits the Pl\"ucker embedding into the projective space $\bP(\Lambda^{a_i})(M_i)$. 
Let $\eO_i(m)$ be the line bundle on $\Gr_{a_i}(M_i)$ obtained as a pullback of $\eO(m)$
on $\bP(\Lambda^{a_i}(M_i))$. 
For $\bfm\in \bZ_{\ge 0}^{Q_0}$ we denote by $\eO(\bfm)$ the exterior tensor product 
$\prod_{i\in Q_0} \eO_i(m_i)$, which is a line bundle on the product of Grassmannians. We use the 
same symbol  $\eO(\bfm)$ to denote the restriction to the quiver Grassmannian $\Gr_{\dim P}(M)$.

Recall the multi-homogeneous ideals ${\mathcal I}(M)$ inside the polynomial ring in Pl\"ucker variables
$\Delta^{(i)}_J$ which defines the reduced scheme structure of the quiver Grassmannians 
$\Gr_{\dim P}(M)$.  The multi-graded Pl\"ucker algebra admits the decomposition
\[
{\rm Pl}(M) = \bC[\Delta_J^{(i)}]/{\mathcal I}(M) = 
\bigoplus_{\bfm\in{\mathbb Z}_{\ge 0}^{Q_0}} {\rm Pl}_\bfm(M),
\]
where the homogeneous component ${\rm Pl}_\bfm(M)$ is spanned by all the monomials of the form
\[
\prod_{i\in Q_0} \Delta^{(i)}_{J(1)}\dots \Delta^{(i)}_{J(m_i)},\ J(r)\subset [d_i],\ 
|J(r)|=a_i\ (r\in [m_i]).
\]

For example, for $Q$ being the equioriented type $A$ quiver the homogeneous components 
${\rm Pl}_\bfm(M^0)$ 
are identified with the dual irreducible highest weight ${\mathfrak{sl}}_n$ modules and  ${\rm Pl}_\bfm(M^1)$
are identified with dual PBW degenerate representations \cite{Fe23}. We put forward the following conjectures.

\begin{conj}\label{conj:dimpl}
Let $P=\bigoplus_{i\in Q_0}$, $I=\bigoplus_{j\in Q^0} I_j$. 	
For any $M\in {\rm Rep}_\bd$ and any $\bfm\in\bZ_{\ge 0}^{Q_0}$ one has the inequality 
$\dim {\rm Pl}_\bfm(M)\ge \dim {\rm Pl}_\bfm(M^0)$. The equality 
$\dim {\rm Pl}_\bfm(M) = \dim {\rm Pl}_\bfm(M^0)$ holds for all $\bfm$ if and only if 
$\dim \Gr_{\dim P}(M) = \langle \dim P,\dim I\rangle$.
\end{conj}

\begin{rem}\label{rem:difdim}
Conjecture \ref{conj:dimpl} does not hold for general $P$ and $I$ even for $Q$ of type $A$. 
More precisely, if one allows arbitrary multiplicities of indecomposable projecive
and injective modules in the decomposition of $P$ and $I$, then the dimensions 
$\dim {\rm Pl}_\bfm(M)$ may be different from $\dim {\rm Pl}_\bfm(M^0)$ even for $M$ degenerating
to $M^1$.
\end{rem}

\begin{conj} \label{conj:hcoh}
Let $P=\bigoplus_{i\in Q_0} P_i$, $I=\bigoplus_{j\in Q^0} I_j$. Then for any $M$ such that 
$\Gr_{\dim P}(M)$ is of minimal dimension, one has for all $\bfm\in\bZ_{\ge 0}^{Q_0}$  
\[
\dim H^0(\Gr_{\dim P},\eO(\bfm)) = \dim {\rm Pl}_\bfm(M),\quad 
H^{>0}(\Gr_{\dim P},\eO(\bfm)) = 0.
\]	
\end{conj}

\begin{rem}
Conjecture  \ref{conj:hcoh} is not expected to hold for general $P$ and $I$ (although we 
do not have a program computing the dimensions of the cohomology groups). The reason comes
from Remark \ref{rem:difdim}: over a flat family the Euler characteristic of $\eO(\bfm)$
is preserved, however the dimensions of the homogeneous components of the Pl\"ucker algebra
vary. This might be an indicator that the higher cohomology groups do not vanish in general. 
\end{rem}

\begin{example}
Let $Q$ be of the form $1\rightarrow 2\leftarrow 3$. Then for $u_i=m_i +1$:
\begin{multline*}
	\dim {\rm Pl}_\bfm(M^0) = 	
\frac{1}{12} u_1u_2u_3(3u_1u_2 + 3u_1u_3 + 3u_2u_3 + 2u_2^2 + 1).
\end{multline*}
\end{example}

\begin{example}
Let $Q$ be of the form 

\begin{picture}(100,110)     
	\put(92,50){\texttt{1}}    
	\put(100,54){\vector(11,0){30}} 
	\put(134,50){\texttt{2}}    
	\put(172,88){\vector(-1,-1){30}} 
	\put(172,18){\vector(-1,1){30}} 
	\put(174,90){\texttt{3}}    
	\put(174,12){\texttt{4}}    
\end{picture}

\noindent Then
\begin{multline*}
\dim {\rm Pl}_\bfm(M^0) = 	
\frac{1}{24}  u_1u_2u_3u_4(3u_1u_2u_3 + 3u_1u_2u_4 + 3u_1u_3u_4 + 3u_2u_3u_4 + \\
u_2^3 + 2u_2^2(u_1 + u_3 + u_4) + 2u_2 + u_1 + u_3 + u_4)
\end{multline*}
(as above, $u_i=m_i+1$).
\end{example}

\section{Software}\label{sec:comp}
Suggesting the above conjectures was made possible after extensive computer experiments with the \texttt{quiver-representation} library we've developed. It's open source and available on github \cite{Fed25}.

The library's main language is Python~3. It also relies on: 

\begin{itemize}
    \item Macaulay 2 \href{https://macaulay2.com/}{https://macaulay2.com/} to derive algebro-geometric properties of quiver Grassmannians from their equations.
    \item GNU Parallel for batch computations.
    \item Graphviz \href{https://graphviz.org/}{https://graphviz.org/} for visualization of $\Gamma_\bd$.
\end{itemize}

The library allows to perform a wide range of module-theoretic computations for quiver representations over the field of complex numbers or over any finite field, including: 

\begin{itemize}
\item Inferring indecomposable simple, projective and injective representations.
\item Computing direct sums, kernels and cokernels of morphisms, radicals and socles, projective covers and injective hulls. 
\item Finding a basis of ${\rm Hom}_Q(M, N)$ for given $M$ and $N$.
\item Enumerating indecomposables of $A_n$ and $D_n$; creating a representation for a bag of intervals.
\end{itemize}

\begin{figure}[htbp]
	\centering
	\includegraphics[width=0.9\linewidth]{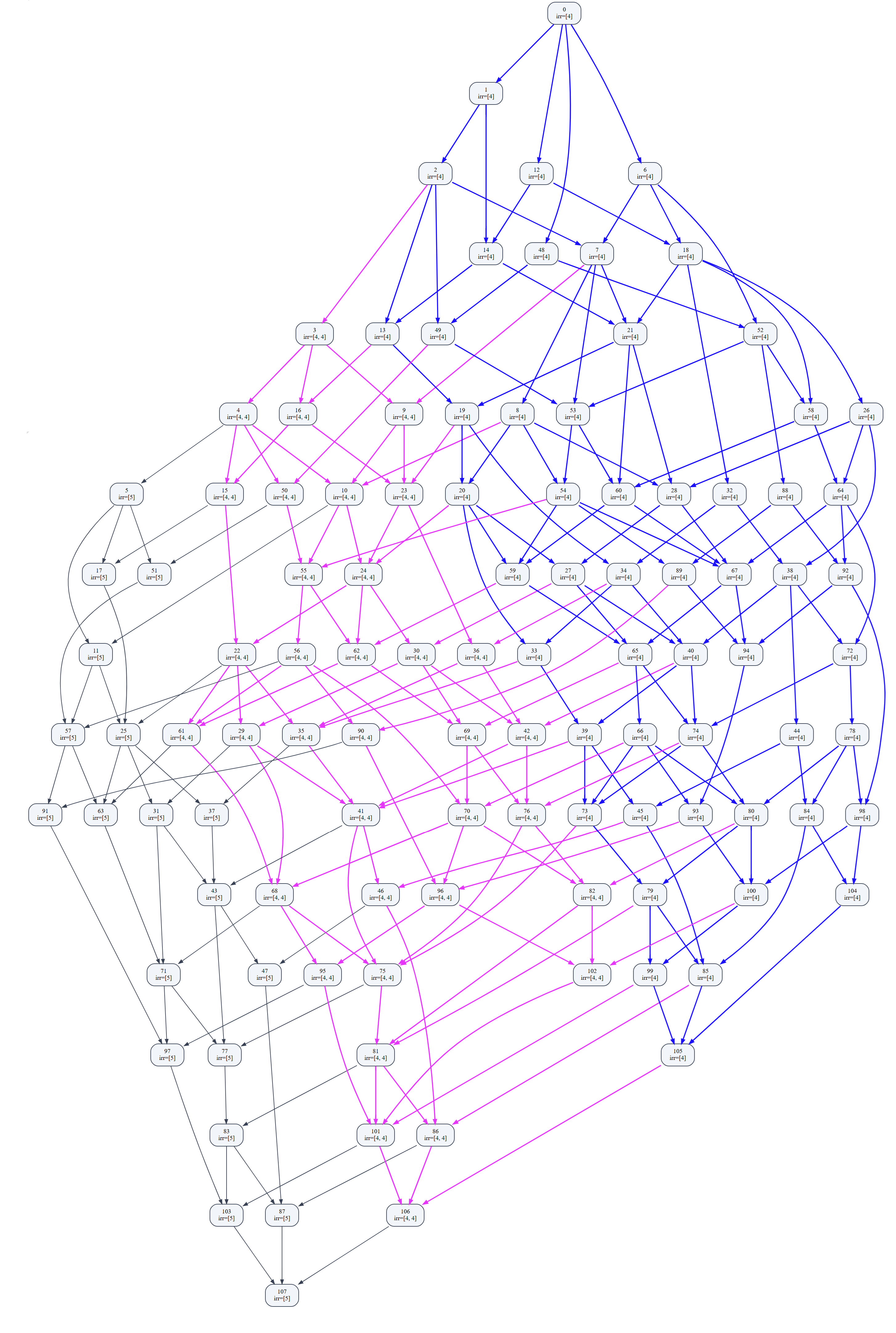}
	\caption{$\Gamma_\bd(2)$ for $Q=\bullet \rightarrow \bullet\leftarrow \bullet\leftarrow \bullet$,
		$P=P_1\oplus P_2\oplus P_3$, $I=I_1\oplus I_2\oplus I_3$.}
\end{figure}

The library also has the full pipeline of checking the hypotheses mentioned in this paper: from enumerating all Pl{\" u}cker and incidence relations for a given quiver Grassmannian to constructing the graph $\Gamma_\bd$ and to computing dimensions of irreducible components and Hilbert functions with Macaulay~2.

Despite the library supports module-theoretic computations over $\mathbb{C}$, we advise you to use it with caution. Floating-point computations are unavoidably imprecise; for example, a tiny computational error might make any matrix full-rank, which would lead to incorrect kernels, cokernels, etc.

Computations over finite fields lack this problem. So, in our experiments, the calculation of $\dim {\rm Hom}_Q(\cdot, \cdot)$ was done over $\mathbb{F}_{107}$.
The data needed to assess our hypotheses using the pipeline can be found in the \texttt{examples} folder: \href{https://github.com/st-fedotov/quiver/tree/main/examples}{https://github.com/st-fedotov/quiver/tree/main/examples}. 

The pipeline tends to suffer from memory exhaustion due to combinatorial explosion when a quiver contains paths of length $\geqslant 5$ or if the multiplicities of indecomposable projectives or injectives in $P\oplus I$ are greater than 1. This can be somewhat relieved by increasing the \texttt{gc\_heap\_size} parameter in the experimental config provided you have enough RAM. Still, the computations might take some time.

The bulk of our experiments were done on a Nebius' virtual machine with 128 CPU and 512Gb RAM. Depending on the memory requirements, we used between 30 and 120 parallel computational streams.
Please check the library's readme for installation guidance and further details.

The library was created with the help of coding agents: GPT-5 in chat mode \href{https://chatgpt.com/}{https://chatgpt.com/}, Codex \href{https://chatgpt.com/codex/}{https://chatgpt.com/codex/}, and Claude Code \href{https://www.claude.com/product/claude-code}{https://www.claude.com/product/claude-code}.

\end{document}